\documentclass[11pt]{amsart}

\usepackage{amsxtra}

\newtheorem{theorem}{Theorem}[section]
\newtheorem{proposition}[theorem]{Proposition}
\newtheorem{lemma}[theorem]{Lemma}

\theoremstyle{definition}

\newtheorem{example}[theorem]{Example}
\theoremstyle{remark}
\newtheorem{remark}[theorem]{Remark}

\numberwithin{equation}{section}

\newcommand{\C}{\mathbb{C}}
\newcommand{\N}{\mathbb{N}}
\newcommand{\R}{\mathbb{R}}

\newcommand{\T}{\mathbb{T}}
\newcommand{\Z}{\mathbb{Z}}
\renewcommand{\phi}{\varphi}
\DeclareMathOperator{\im}{Im}
\DeclareMathOperator{\re}{Re}

\begin{document}

\title[Periodically modulated leaky wires]{Absolute continuity of the
  spectrum for periodically modulated leaky wires in $\R^3$}

\author{Pavel Exner}
\address{Department of Theoretical Physics, Nuclear Physics
Institute, Academy of Sciences, 25058 Rez near Prague, Czech
Republic}
\email{exner@ujf.cas.cz}

\author{Rupert L. Frank}
\address{Royal Institute of Technology, Department of Mathematics,
  Lindstedts\-v\"agen 25, 10044 Stockholm, Sweden}
\email{rupert@math.kth.se}

\begin{abstract}
We consider a model of leaky quantum wire in three dimensions. The
Hamiltonian is a singular perturbation of the Laplacian supported
by a line with the coupling which is bounded and periodically
modulated along the line. We demonstrate that such a system has a
purely absolutely continuous spectrum and its negative part has
band structure with an at most finite number of gaps. This result
is extended also to the situation when there is an infinite number
of the lines supporting the perturbations arranged periodically in
one direction.
\end{abstract}

%\keywords{Scattering theory, periodic operator, Schr\"odinger
%operator, singular potential}
%\subjclass[2000]{Primary 35J10; Secondary 35J25, 35P05, 35P25}

\date{\today}

\maketitle

%%%%%%%%%%%%%%%%%%%%%%%%%%%%%%%%%%%%%%%%%%%%%%%%%%%%%%%%%%%%%%%%%%%%%%%%%%%%

\section{Introduction}

Existence of transport in quantum systems having a periodic
structure is important in many areas, particularly in
condensed-matter physics. Mathematically this property is
expressed as the absolute continuity of spectrum of the
appropriate Hamiltonian. For Schr\"odinger operators with regular
potentials which are ``completely'' periodic in the sense that the
basic period cell is compact the problem is well understood --
cf.~\cite{T} or \cite{RS}, Sec.XIII.16.

Recently a class of models attracted attention in which one or
both of the above conditions are violated. They concern thin
microscopic semiconductor structures, often dubbed ``quantum
wires'', which are intensively studied as construction elements of
future electronic devices. Comparing to the usual treatment of
such objects the mentioned models are realistic in the sense that
they describe the wires by elongated ``potential wells'' so that
quantum tunneling is not suppressed. On the other hand, they are
often idealized using singular potentials with the aim to make the
model solvable; the corresponding Hamiltonian can be written
formally as
\begin{equation}\label{eq:formal}
-\Delta + \sigma(x) \delta(x-\Gamma)
\end{equation}
in $L_2(\R^\nu),\; \nu=2,3$, where
$\Gamma$, typically a curve or a family of curves, supports the
interaction. Various spectral properties of such operators have
been derived in several last years, see \cite{AGHH} for a
bibliography.

On the other hand, many questions are still open. For instance,
while it is natural to conjecture that the spectrum is absolutely
continuous when $\Gamma$ is a periodic curve and $\sigma$ is
constant along it, only a partial result is known \cite{BeDuE} and
a full proof is missing. The situation is better in the case when
the periodicity concerns the coupling rather than the geometry of
the interaction support. If $\Gamma$ is a straight line in $\R^2$
and $\sigma$ is periodic the sought property can be obtained by
modification of the results of \cite{Fr1, Fr2, FrSh} (recall also
that a similar result for Schr\"odinger operators with a
``partially periodic'' regular potential was derived in
\cite{FiKl1}). It is important, however, that the codimension of
$\Gamma$ is one here and the Hamiltonian can be defined naturally
through the associated quadratic form.

The main aim of this paper is to solve the analogous problem for a
straight line in $\R^3$. In this case the codimension is two and
the operator has to be defined by means of boundary conditions
involving generalized boundary values as in \cite{EKo1}. We will be
able to demonstrate that such a Hamiltonian has a purely
absolutely continuous spectrum, and moreover, that its negative
part has band structure with an at most finite number of gaps. The
corresponding generalized eigenfunctions are of physical interest,
of course, because they describe states guided along the ``wire''.
Moreover, we are going to extend the result to the case where the
interaction support consists of an infinite family of parallel and
equidistant straight lines in a fixed plane\footnote{We leave out
in this paper another possible extension to the situation when the
line family is periodic in two different directions. In this case
the basic period cell is compact and more conventional methods can
be used.}.

The paper is organized as follows. The main results are formulated
in the next section and proved subsequently in
Secs.~\ref{sec:defop}--\ref{sec:pointspec}. The results concerning
the extended model are stated in Sec.~\ref{sec:model2} and proved
in the rest of the paper.

%%%%%%%%%%%%%%%%%%%%%%%%%%%%%%%%%%%%%%%%%%%%%%%%%%%%%%%%%%%%%%%%%%%%%%%%%%%%
\
\section{Main results}\label{sec:mainres}

\subsection{Description of the results}\label{sec:mainres1}
In $\R^3$ we introduce coordinates $(x,y)$, $x\in\R,\ y\in\R^2$,
and denote $\Gamma:= \R\times\{(0,0)\}$. Moreover, let $\sigma$ be
a real-valued, $2\pi$-periodic function such that
\begin{equation}\label{eq:sigma}
    \sigma\in L_\infty(\R).
\end{equation}
We will construct a  self-adjoint operator $H$ in $L_2(\R^3)$
corresponding to the formal expression (\ref{eq:formal}) which can
be written\footnote{It has to be stressed that this expression is
formal and the proper way to introduce $\sigma$ is given by
(\ref{eq:defoppi}) below; recall that the \emph{absence of the
coupling means} $\sigma=\infty$.} also as
$-\Delta+\sigma(x)\delta(|y|)$. Put $C_\epsilon:=\{(x,y)\in\R^3:\
|y|\leq\epsilon\}$ for $\epsilon>0$. We consider functions $u\in
L_2(\R^3)$ such that
\begin{equation}\label{eq:h2}
    u\in H^2(\R^3\setminus C_\epsilon) \qquad\text{for all
    }\epsilon>0.
\end{equation}
By the embedding theorems $u$ is continuous in $\R^3\setminus\Gamma$
and hence it restriction $u(.,y)$ to the line $\{(x,y):\,x\in\R\}$ is
well-defined. We denote by $\Upsilon$ the class of functions $u\in
L_2(\R^3)$ satisfying \eqref{eq:h2} and such that the limits
\[ \Xi u := -\lim_{y\to 0}
 \frac1{\log |y|}\ u(.,y), \qquad
 \Omega u := \lim_{y\to 0}
 \left(u(.,y) +\log |y| \Xi u\right), \]
exist in the sense of distributions and belong to $L_2(\R)$.\\
We will recall in Subsection \ref{sec:defop1} how one constructs a self-adjoint
operator $H$ in $L_2(\R^3)$ such that
\begin{equation}\label{eq:defop}\begin{split}
    H u & = -\Delta u \qquad\text{in }\R^3\setminus\Gamma,\\
    \mathcal{D}(H) & =\{u\in\Upsilon: \Delta u\in
        L_2(\R^3\setminus\Gamma),\ \Omega u - 2\pi\sigma\Xi u = 0 \}.
    \end{split}
\end{equation}
(Recall that $\Delta u \in L_2(\R^3\setminus\Gamma)$ means that
the distribution $\Delta u$ is a function on $\R^3\setminus\Gamma$
and square-integrable. We do not make an assertion about its
nature on $\Gamma$.)\\
Our main result is

\begin{theorem}\label{thmac} The spectrum of the operator $H$ is
purely absolutely continuous.
\end{theorem}

For the proof of this theorem we will investigate the scattering
between $H$ and $H_0$, the standard self-adjoint realization of $-\Delta$ in $\R^3$.
Recall the definition (in case of existence) of the \emph{wave operators} (see, e.g.,
\cite{Ya})
\begin{equation}\label{eq:defwo}
W_\pm := s-\lim_{t\rightarrow\pm\infty}\exp(itH)\exp(-itH_0).
\end{equation}
We will prove

\begin{theorem}\label{thmscat}
The wave operators $W_+$ and $W_-$ exist, satisfy
$\mathcal{R}(W_+)=\mathcal{R}(W_-)$ and are not complete.
\end{theorem}

The existence of the wave operators implies, as it is well-known,
that $\sigma_{ac}(H_0)=[0,\infty)$ is contained in the absolutely
continuous spectrum of the operator $H$. Moreover, we note that
the identity $\mathcal{R}(W_+)=\mathcal{R}(W_-)$ implies the
unitarity of the scattering matrix. The non-completeness of the
wave operators is due to \emph{guided states}, i.e. states that
are localized near the wire $\Gamma$ for all times. They
correspond to bands in the (negative) spectrum of $H$. We will
prove

\begin{theorem}\label{thmsurf}
The negative spectrum of the operator $H$ is non-empty and has
band structure with at most finitely many gaps.
\end{theorem}

However, we emphasize that guided states correspond not only to
negative energies. Indeed, if $\sigma\equiv\alpha\in\R$ is
constant then the spectrum of $H$ on $\mathcal{R}(W_\pm)^\bot$
coincides with the half-line $[\xi(\alpha),\infty)$ where
\begin{equation}\label{eq:xi}
    \xi(\alpha)= -4e^{2(-2\pi\alpha+\psi(1))}, \qquad\alpha\in\R,
\end{equation}
and $-\psi(1)$ is the Euler constant (numerically,
$-\psi(1)=0.577...$). It is a natural question whether the
spectrum of $H$ on  $\mathcal{R}(W_+)^\bot$ is bounded when
$\sigma$ is non-constant\footnote{In the case of a constant
$\sigma$ the positive-energy guided states are expected be
unstable with respect to perturbations, but we are not going to
discuss this problem here.}.

\begin{remark} The subspace $\mathfrak{M}$ of functions being rotationally
symmetric with respect to the variable $y$ reduces both $H$ and $H_0$,
and the parts of these operators in $\mathfrak{M}^\bot$ coincide. In
particular, guided states belong to $\mathfrak{M}$.
\end{remark}

\begin{remark} Of course, the assumption that the period is $2\pi$
is not essential. Moreover, the assumption \eqref{eq:sigma} can be
relaxed to assuming that for every $\epsilon>0$ there is a
$C_\epsilon>0$ such that
\[
\|\sigma f\|^2_{L_2(-\pi,\pi)} \leq \epsilon \sum_{n\in\Z}
\left(\log(1+n^2)\right)^2 |\hat f_n|^2 + C_\epsilon
\|f\|^2_{L_2(-\pi,\pi)}
\]
for all smooth, $2\pi$-periodic functions $f$ with Fourier
coefficients $\hat f_n$, see \eqref{eq:fourier} below.
\end{remark}

\subsection{Direct integral decomposition}\label{sec:directintegral}
Because of periodicity the operator $H$ can be partially
diagonalized. A fundamental cell is the layer
\[ \Pi := \{ (x,y)\in\R^3:\; x\in [-\pi,\pi)\}. \]
Actually, we will only work with functions on $\Pi$ that are periodic
with respect to the variable $x$ and one may think of $\Pi$ as a manifold
with opposite points on the planes $\{x=\pi\}$ and $\{x=-\pi\}$ identified,
but we ignore this for the sake of simplicity. However, we will identify
$\Pi\cap\Gamma= [-\pi,\pi)\times\{(0,0)\}$ with the 'torus' $\T$.\\
By $\tilde H^2(\Pi)$ we denote the class of functions $u\in
H^2(\Pi)$ the periodic extension of which belongs to
$H^2_{loc}(\R^3)$, and by $\tilde\Upsilon$ we denote the class
of functions $u\in L_2(\Pi)$ satisfying
\[
u\in H^2(\Pi\setminus C_\epsilon)
\qquad\text{for all }\epsilon>0
\]
and such that their periodic extension belongs to $\Upsilon_{loc}$. Here as
usual if $F$ is a class of functions on $\R^3$ then
$F_{loc}:=\{u:\R^3\rightarrow\C : \phi u\in F\ \forall \phi\in
C_0^\infty(\R^3) \}$. For $u\in\tilde\Upsilon$ the functions $\Xi u$,
$\Omega u$ are well-defined and belong to $L_2(\T)$.\\
We will recall in Subsection \ref{sec:defoppi} that there exists a
family of self-adjoint operators $H(k),\ k\in
Q:=[-\frac{1}{2},\frac{1}{2})$, in $L_2(\Pi)$ such that
\begin{equation}\label{eq:defoppi}\begin{split}
    H(k) u & =((D_x+k)^2+D_y^2) u
    \qquad\text{in }\Pi\setminus\Gamma,\\
    \mathcal{D}(H(k)) & =\{u\in\tilde\Upsilon:\, \Delta u\in
    L_2(\Pi\setminus\Gamma),\ \Omega u - 2\pi\sigma\Xi u = 0 \}.
\end{split}\end{equation}
Here $D_x=-i\frac{\partial}{\partial x},\ D_y = -i\nabla_y$.
Moreover, we denote by $H_0(k)$ the operator $(D_x+k)^2+D_y^2$ in $L_2(\Pi)$
with domain $\tilde H^2(\Pi)$.\\
The \emph{Gelfand transformation} is initially defined for $u\in
C^{\infty}_0(\R^3)$ by
\[ (\mathcal{U}u)(k,x,y) := \sum_{n\in\Z} e^{-i k(x+2\pi n)} u(x+2\pi n,y),
        \qquad k\in{\textstyle Q},\, (x,y)\in\Pi, \]
and extended by continuity to a \emph{unitary} operator
$\mathcal{U}:L_2(\R^3) \rightarrow \int_Q^\oplus L_2(\Pi)\,dk$. It is
well-known that
\begin{equation}\label{eq:directintdecompunperturbed}
    \mathcal{U}\, H_0\, \mathcal{U}^* = \int_Q^\oplus H_0(k)\ dk.
\end{equation}
In Subsection \ref{sec:defoppi} we will prove that similarly
\begin{equation}\label{eq:directintdecomp}
    \mathcal{U}\, H\, \mathcal{U}^* = \int_Q^\oplus H(k)\ dk.
\end{equation}
This reduces the investigation of the operator $H$ to the study of the
fiber operators $H(k)$.

\subsection{Results about the fiber operators}\label{sec:scatpi}
Information about the continuous spectrum of the operators $H(k)$ can
be obtained by scattering theory for the pair $(H(k),H_0(k))$.
Note that the operator $H_0(k)$ can be diagonalized explicitly.
Its spectrum is purely absolutely continuous and coincides with
$[k^2,\infty)$. The spectral multiplicity is finite and changes at
the points from the \emph{threshold set}
\[ \tau(k) := \{ (n+k)^2: n\in\Z\}. \]
We introduce the wave operators
\[ W_\pm(k) := s-\lim_{t\rightarrow\pm\infty}\exp(itH(k))\exp(-itH_0(k)). \]

\begin{proposition}\label{thmacpi} Let $k\in Q$. Then the wave
operators $W_+(k)$ and $W_-(k)$ exist and are complete. In particular,
$\sigma_{ac}\left(H(k)\right) = [k^2,\infty)$.
\end{proposition}

By deriving a limiting absorption principle we will show

\begin{proposition}\label{thmscpi} Let $k\in Q$. Then $\sigma_{sc}\left(H(k)\right) =\emptyset$.
\end{proposition}

Concerning the point spectrum of the fiber operators we prove

\begin{proposition}\label{thmdiscpi} Let $k\in Q$. Then $\sigma_{disc}\left(H(k)\right)
=\sigma\left(H(k)\right)\cap (-\infty,k^2)$ is non-empty and
finite.
\end{proposition}

Indeed, we will prove that $H(k)$ has an eigenvalue less or equal
$\xi(\tilde\sigma)+k^2$ where $\xi$ is given by \eqref{eq:xi} and
\begin{equation}\label{eq:sigmaaverage}
    \tilde\sigma := \frac{1}{2\pi}\int_\T \sigma(x)\ dx.
\end{equation}
Moreover, for the proof of absolute continuity we need

\begin{proposition}\label{thmpppi} There exists a countable family of
open connected sets $U_j,\ V_j\subset\R$ and real-analytic
functions $h_j:U_j\times V_j\rightarrow \C$ satisfying
\begin{enumerate}
\item for all $j$ and all $\lambda\in U_j$ one has $h_j(\lambda,.)\not\equiv
0$, and
\item $\{ (\lambda,k)\in\R\times Q :
\lambda\in\sigma_{p}\left(H(k)\right)\} \subset \bigcup_j \{ (\lambda,k)\in U_j\times
V_j : h_j(\lambda,k)=0 \}$.
\end{enumerate}
\end{proposition}

\subsection{Reduction to the fiber operators}\label{sec:reduction}
Our main results Theorems \ref{thmac} - \ref{thmsurf} can be
deduced from Propositions \ref{thmacpi} - \ref{thmpppi} in a
standard way. We only sketch the major steps.

\begin{proof}[Proof of Theorem \ref{thmac}] Propositions
\ref{thmacpi}, \ref{thmscpi} and \ref{thmpppi} allow us to follow
the proof of Theorem 1.4 in \cite{Fr2} word by word.
\end{proof}

\begin{proof}[Proof of Theorem \ref{thmsurf}] We will see below
that the discrete eigenvalues of $H(k)$ depend piecewise
analytically on $k$. Hence the existence of negative spectrum of
$H$ and its band structure follow from Proposition \ref{thmdiscpi}
and the decomposition \eqref{eq:directintdecomp}. By analytic
perturbation theory there can be at most one gap in
$\sigma\left(H\right)$ between two consecutive eigenvalues of
$H(0)$. Hence the finiteness of gaps follows again from
Proposition \ref{thmdiscpi}.
\end{proof}

\begin{proof}[Proof of Theorem \ref{thmscat}] Proposition
\ref{thmacpi} implies the existence of the wave operators $W_\pm$
and the equality of their ranges (see \cite{Fr1}). The
non-complete\-ness follows immediately from Theorems \ref{thmac} and
\ref{thmsurf}.
\end{proof}

%%%%%%%%%%%%%%%%%%%%%%%%%%%%%%%%%%%%%%%%%%%%%%%%%%%%%%%%%%%%%%%%%%%%%%%%%%%%%%

\section{Definition of the operators}\label{sec:defop}

\subsection{Definition of the operator $H$}\label{sec:defop1}
Recall that the domain of $H_0$ is $H^2(\R^3)$ and that the trace
operator $\gamma:H^2(\R^3)\rightarrow L_2(\R)$,
\[ \gamma u := u|_\Gamma, \]
is well-defined. (Here we identify $\Gamma$ naturally with
$\R$.)\\
For $z\in\C\setminus [0,\infty)$ we consider the pseudo-differential
  operator $T(z)$ in $L_2(\R)$,
\[\begin{split}
T(z) & := \frac{1}{4\pi}\log\left(D^2-z\right) -\varsigma I,\\
\mathcal{D}(T(z)) & := \{f\in L_2(\R):\  \int_\R
\left(\log(1+\xi^2)\right)^2 |\hat{f}(\xi)|^2\ d\xi <\infty \},
\end{split}\]
where $\varsigma=\frac{1}{2\pi}(\ln 2+\psi(1))$ and $-\psi(1)$ is
as before the Euler constant. Here and in all the following we
choose the principal branch of the logarithm on
$\C\setminus(-\infty,0]$. (Note that we have changed the sign in
the
definition of $T(z)$ as compared to \cite{EKo1}.)\\
We write $R_0(z):=(H_0-zI)^{-1}$. One checks easily that for
$z,\zeta\in\C\setminus[0,\infty)$
\begin{equation}\label{eq:asspo}
    T(\overline z)=T(z)^*, \qquad
    T(z)-T(\zeta)= (\zeta-z)\left(\gamma R_0(\zeta)\right)
    \left(\gamma R_0(\overline z)\right)^*.
\end{equation}
By abstract arguments of \cite{Po} (see also \cite{EKo1}) this implies
that there exists a self-adjoint operator $H$ in $L_2(\R^3)$ such that
\[
\{ z\in\C\setminus[0,\infty):\ 0\in\rho\left(T(z)+\sigma\right) \}
\subset\rho(H)
\]
and such that the resolvent $R(z):=(H-zI)^{-1}$ is related to $R_0(z)$
by
\begin{equation}\label{eq:resdiff}\begin{split}
    R(z)= & R_0(z) + (\gamma R_0(\overline z))^*(T(z)+\sigma)^{-1}
    \gamma R_0(z),\\
    & z\in\C\setminus[0,\infty),\ 0\in\rho\left(T(z)+\sigma\right).
\end{split}\end{equation}
By \eqref{eq:sigma} the operator $T(-a)+\sigma$ is positive definite for all
sufficiently large $a$, and hence $H$ is lower semibounded. Moreover,
it was shown in \cite{EKo1} that the operator defined in this way
satisfies \eqref{eq:defop}. Without reproducing the proof here we note
that it relies on the identities
\[
\Xi \left(\gamma R_0(\overline z)\right)^* = \frac{1}{2\pi}I,
\qquad \Omega \left(\gamma R_0(\overline z)\right)^* = - T(z).
\]

\subsection{Definition of the operators $H(k)$}\label{sec:defoppi}
Recall that functions in the domain $\tilde H^2(\Pi)$ of $H_0(k)$
satisfy periodic boundary conditions with respect to the variable $x$,
and that we identify $\Gamma\cap\Pi$ with $\T$. We use the same
notation $\gamma$ for the trace operator $\tilde H^2(\Pi)\rightarrow
L_2(\T)$.\\
Fix $k\in Q$. For $z\in\C\setminus [k^2,\infty)$ we consider the
pseudo-differential operator $T(z,k)$ in $L_2(\T)$,
\[\begin{split}
T(z,k) & := \frac1{4\pi}\log\left((D+k)^2-z\right) -\varsigma I,\\
\mathcal{D}(T(z,k)) & := \{f\in L_2(\T):\ \sum_{n\in\Z}
\left(\log(1+n^2)\right)^2 |\hat{f}_n|^2 <\infty \}.
\end{split}\]
Here
\begin{equation}\label{eq:fourier}
    \hat f_n:=\frac{1}{\sqrt{2\pi}}\int_\T f(x) e^{-inx} dx,
    \qquad n\in\Z,
\end{equation}
denote the Fourier coefficients of $f\in L_2(\T)$, in terms of which
the action of $T(z,k)$ is given by
\[
\left(T(z,k)f\right)\sphat_n =
\left(\frac1{4\pi}\log\left((n+k)^2-z\right)-\varsigma\right) \hat
f_n,
\qquad n\in\Z.
\]
Our next goal is to construct the operators $H(k)$ and to verify the
direct integral decomposition \eqref{eq:directintdecomp}. For this we
introduce the unitary operator $\tilde{\mathcal
  U}:L_2(\R)\to\int_Q^\oplus L_2(\T)\,dk$, defined for $f\in
C_0^\infty(\R)$ by
\[
(\tilde{\mathcal U}f)(k,x) := \sum_{n\in\Z} e^{-i k(x+2\pi n)} f(x+2\pi n),
    \qquad k\in{\textstyle Q},\, x\in\T.
\]
We note that on $H^2(\R^3)$ one has the identity
\begin{equation*}\label{eq:gelfandtrace}
  \gamma\,\mathcal U = \tilde{\mathcal U}\,\gamma,
\end{equation*}
where $\gamma$ is the trace operator $\int_Q^\oplus \tilde
H^2(\Pi)\,dk\rightarrow \int_Q^\oplus L_2(\T)\,dk$ on the LHS and the
trace operator  $H^2(\R^3)\rightarrow L_2(\R)$ on the RHS. We denote
the 'unperturbed' resolvent by $R_0(z,k):=(H_0(k)-zI)^{-1}$. In view
of \eqref{eq:directintdecompunperturbed} we find
\begin{equation}\label{eq:gelfandfree}
    \tilde{\mathcal U}\gamma R_0(z)\,\mathcal U^* =
    \int_Q^\oplus \gamma R_0(z,k)\,dk,
    \qquad z\in\C\setminus[0,\infty).
\end{equation}
Moreover, it turns out that
\begin{equation}\label{eq:gelfandt}
  \tilde{\mathcal U}\, T(z)\, \tilde{\mathcal U}^* =
  \int_Q^\oplus T(z,k)\,dk, \qquad z\in\C\setminus[0,\infty).
\end{equation}
Combining \eqref{eq:gelfandfree}, \eqref{eq:gelfandt} with
\eqref{eq:asspo} we conclude easily that
\begin{equation*}
    T(\overline z,k)=T(z,k)^*, \qquad
    T(z,k)-T(\zeta,k)= (\zeta-z)\left(\gamma
    R_0(\zeta,k)\right)\left(\gamma R_0(\overline z,k)\right)^*
\end{equation*}
for all $k\in Q$, $z,\zeta\in\C\setminus[k^2,\infty)$. (Originally,
  these relations follow for $z\in\C\setminus[k^2,\infty)$ and can be
  extended by analyticity to $z\in[0,k^2)$. Alternatively, they
  may be established directly.) It follows again from \cite{Po} that there
  exists a self-adjoint operator $H(k)$ in $L_2(\Pi)$ such that its
  resolvent $R(z,k):=(H(k)-zI)^{-1}$ satisfies
\begin{equation}\label{eq:resdiffpi}\begin{split}
    R(z,k)= & R_0(z,k)+(\gamma R_0(\overline
    z,k))^*(T(z,k)+\sigma)^{-1}\gamma R_0(z,k),\\
    & z\in\C\setminus[k^2,\infty),\
    0\in\rho\left(T(z,k)+\sigma\right).
    \end{split}
\end{equation}
Combining this with \eqref{eq:directintdecompunperturbed},
\eqref{eq:gelfandfree}, \eqref{eq:gelfandt} we obtain
\[
\mathcal U R(z)\,\mathcal U^* = \int_Q^\oplus R(z,k)\,dk,
\qquad z\in\rho(H),
\]
which implies \eqref{eq:directintdecomp}. Finally, the
characterization \eqref{eq:defoppi} is deduced from \eqref{eq:defop}
as in \cite{EKo2}.

%%%%%%%%%%%%%%%%%%%%%%%%%%%%%%%%%%%%%%%%%%%%%%%%%%%%%%%%%%%%%%%%%%%%%%%%%%%%%%

\section{The continuous spectrum of the operators $H(k)$}\label{sec:contspec}

\subsection{Proof of Proposition \ref{thmacpi}}\label{sec:thmacpi}
According to a result of Birman-Kre\u\i n (see \cite{Ya}) it suffices to prove that
\[ R(z_0,k)-R_0(z_0,k)\in\mathfrak{S}_1, \]
the trace class, for some
$z_0\in\rho\left(H(k)\right)\cap\rho\left(H_0(k)\right)$. For
$a>0$ sufficiently large the operator $T(-a,k)+\sigma$ is boundedly
invertible. In view of \eqref{eq:resdiffpi} it suffices therefore
to prove that
\[ \gamma R_0(-a,k)\in\mathfrak{S}_2, \]
the Hilbert-Schmidt class. For this recall that $R_0(z,k)$,
$z\in\C\setminus[0,\infty)$, is an integral operator with the
kernel
\begin{equation}\label{eq:freekernel}\begin{split}
    & r_0(x,y,x',y',z) := \frac{1}{(2\pi)^2}\sum_{n\in\Z} e^{in(x-x')}\,
        K_0(\sqrt{(n+k)^2-z}\,|y-y'|),\\
    & \qquad (x,y),(x',y')\in\Pi, \end{split}\end{equation}
where $K_0$ is the Macdonald function (or modified Bessel function
of the second kind) of order zero (see \cite{AS}). Here and in the
following we choose the branch of the square root on $\C\setminus
(-\infty,0]$ satisfying $\re\sqrt{.}>0$. It follows using
Parseval's identity that
\[ \begin{split} & \|\gamma R_0(-a,k)\|_{\mathfrak{S}_2}^2 = \\
    & \qquad = \frac{1}{(2\pi)^4} \int_\T \int_\T \int_{\R^2}
        \left|\sum_{n\in\Z}
        e^{in(x-x')}\,K_0(\sqrt{(n+k)^2+a}\,|y|)\right|^2
        dy\, dx\, dx' =\\
    & \qquad =\frac{1}{2\pi} \sum_{n\in\Z} \frac{1}{(n+k)^2+a}
        \int_0^\infty |K_0(r)|^2\, r\,dr.
\end{split} \]
Since the last integral is finite by properties of the Macdonald
function, the proof of Proposition \ref{thmacpi} is complete.

\begin{remark}\label{spclass} A more careful analysis shows that the
  singular values $s_j$ of the operator $R(z_0,k)-R_0(z_0,k)$ satisfy
  the estimate $\sup_j j^{1/p} s_j <\infty$ with $p=\frac{1}{2}$. This
  should becompared with the exponent $p=\frac{2}{3}$ when the
  perturbation of $H_0$ is supported on a two-dimensional plane. (This
  can be seen as in the proof of Corollary 3.3 in \cite{Fr1}.)
\end{remark}

\subsection{The limiting absorption principle for the unperturbed
operator $H_0(k)$}\label{sec:lapfree}
In order to prove Proposition \ref{thmscpi} we have to study the
behaviour of the resolvent $R(z,k)$ as the spectral parameter
approaches the real axis. The relation \eqref{eq:resdiffpi} suggests
that we begin
with the investigation of the unperturbed resolvent $R_0(z,k)$.
Recall the definition of the threshold set $\tau(k)$ where the
spectral multiplicity of the operator $H_0(k)$, and according to
Proposition \ref{thmacpi} also of $H(k)$, changes. Denote by
$\Lambda$ the operator of multiplication by the function
$(1+|y|^2)^{-1/2}$ in $L_2(\Pi)$.

\begin{lemma}\label{lapfree} Let $k\in Q,\
\lambda\in\R\setminus\tau(k)$ and $s>\frac{1}{2}$. Then the
operators
\[ \Lambda^s R_0(\lambda\pm i\epsilon,k) \Lambda^s, \qquad
     \gamma R_0(\lambda\pm i\epsilon,k) \Lambda^s, \qquad \epsilon>0, \]
have limits in Hilbert-Schmidt norm as $\epsilon\rightarrow 0+$.
The convergence is uniform in $\lambda$ from compact subsets of
$\R\setminus\tau(k)$.
\end{lemma}

\begin{proof} This follows in a straightforward way from the
explicit expression \eqref{eq:freekernel} and standard properties
of the Bessel function involved.
\end{proof}

\subsection{Proof of Proposition \ref{thmscpi}}\label{sec:thmscpi}
Let $k\in Q$ be fixed. For $z\in\C_+$ the operator $(T(z,k)-i)^{-1}$
exists, is compact and depends analytically on $z$. Moreover, for
any $\lambda\in\R\setminus\tau(k)$ this family has an analytic
extension to a neighbourhood of $\lambda$ in
$\overline{\C_-}$. Applying the analytic Fredholm alternative (see
Theorem VII.1.9 in \cite{K}) to the family $(T(z,k)-i)^{-1}(\sigma+i)$
we conclude that there exists a set $\mathcal{N}_+(k)$, discrete in
$\R\setminus\tau(k)$, such that the operators
\[ \left(T(\lambda+i\epsilon,k)+\sigma\right)^{-1} =
    \left(I+(T(\lambda+i\epsilon,k)-i)^{-1}(\sigma+i)\right)^{-1}
    (T(\lambda+i\epsilon,k)-i)^{-1}
    \]
have a bounded limit as $\epsilon\rightarrow 0+$ for all
$\lambda\in\R\setminus\left(\tau(k)\cup\mathcal{N}_+(k)\right)$. The
limit is uniform for $\lambda$ from compact subsets of this set.\\
Combining this with relation \eqref{eq:resdiffpi} and Lemma
\ref{lapfree} we see that the operators
\[ \Lambda^s R(\lambda+ i\epsilon,k) \Lambda^s, \qquad \epsilon>0, \]
have limits as $\epsilon\rightarrow 0+$ for all
$\lambda\in\R\setminus\left(\tau(k)\cup\mathcal{N}_+(k)\right)$ (in
Hilbert-Schmidt norm). Moreover, the limit is uniform for $\lambda$
from compact subsets of this set. This implies (see \cite{RS}) that
$\sigma_{sc}\left(H(k)\right)\subset\overline{\tau(k)\cup\mathcal{N}_+(k)}$,
and since the latter set is countable the assertion of Proposition
\ref{thmscpi} follows.

\begin{remark}\label{setn} Denote by $T(\lambda+i0,k)$,
$\lambda\in\R\setminus\tau(k)$,
the boundary value of the operator function $T(z,k)$, $z\in\C_+$.
Then $\mathcal{N}_+(k)$ consists of the values
$\lambda\in\R\setminus\tau(k)$ such that $-1$ is an eigenvalue of
the operator $T(\lambda+i0,k)^{-1}\sigma$. In the next
subsection we will see that this is equivalent to $\lambda$ being
a (non-threshold) eigenvalue of $H(k)$.
\end{remark}

%%%%%%%%%%%%%%%%%%%%%%%%%%%%%%%%%%%%%%%%%%%%%%%%%%%%%%%%%%%%%%%%%%%%%%%%%%%%%%%

\section{The point spectrum of the operators
$H(k)$}\label{sec:pointspec}

Since the perturbed operator in our case is not defined via a quadratic
form we cannot use the usual Birman-Schwinger principle for the
study of the eigenvalues of $H(k)$. In Subsection \ref{sec:opbdry}
we will prove a convenient substitute.

\subsection{Characterization of eigenvalues of $H(k)$}\label{sec:opbdry}
Let $k\in Q$, $\lambda\in\R\setminus\tau(k)$ and define
\begin{equation}\label{eq:alpha}
    \alpha_n(\lambda,k) :=
    \frac{1}{4\pi}\log\left|(n+k)^2-\lambda\right|-\varsigma,
    \qquad n\in\Z.
\end{equation}
In the Hilbert space $L_2(\T)$ we consider the operator
\begin{equation}\label{eq:defopa}\begin{split}
    (A(\lambda,k)f)(x) & := \frac{1}{\sqrt{2\pi}}\sum_{n\in\Z}
    \alpha_n(\lambda,k)\hat f_n e^{inx} +\sigma(x)f(x), \qquad
    x\in\T,\\
    \mathcal{D}(A(\lambda,k)) & := \{f\in L_2(\T): \sum_{n\in\Z}
    \left(\log(1+n^2)\right)^2 |\hat{f}_n|^2 <\infty \}.
    \end{split}
\end{equation}
In the case $\sigma\equiv0$ we will denote this operator by
$A_0(\lambda,k)$. Note that the operator in this case differs from
the operator $T(\lambda+i0,k)$ from Remark~\ref{setn} only on the
subspace $\{f\in L_2(\T) :\ \hat f_n=0,\ (n+k)^2 <\lambda \}$. The
definition on that subspace is rather arbitrary (see Remark
\ref{lowestmodes}) and chosen only for technical convenience.\\
The compactness of the embedding of $\mathcal{D}(A(\lambda,k))$ in
$L_2(\T)$ implies that the operator $A(\lambda,k)$ has compact
resolvent.\\
Now we characterize the non-threshold eigenvalues of the operator
$H(k)$ as the values $\lambda$ for which $0$ is an eigenvalue of
the operator $A(\lambda,k)$. More precisely, we have

\begin{proposition}\label{evchar} Let $k\in Q$ and
  $\lambda\in\R\setminus\tau(k)$.
\begin{enumerate}
\item Let $u\in\mathcal{N}(H(k)-\lambda I)$ and define
\begin{equation}\label{eq:f-u}
        f := \Xi u.
\end{equation}
Then $f\in\mathcal{N}(A(\lambda,k))$, $\hat f_n=0$ if
$(n+k)^2<\lambda$ and, moreover,
\begin{equation}\label{eq:u-f}
        u(x,y) = \frac{1}{\sqrt{2\pi}}\sum_{(n+k)^2>\lambda}
        \hat f_n\,e^{inx}\, K_0(\sqrt{(n+k)^2-\lambda}\ |y|),
        \qquad (x,y)\in\Pi\setminus\Gamma.
\end{equation}
\item Let $f\in\mathcal{N}(A(\lambda,k))$ such that $\hat f_n=0$
if $(n+k)^2<\lambda$ and define $u$ by \eqref{eq:u-f}.\\
Then $u\in\mathcal{N}(H(k)-\lambda I)$ and, moreover,
\eqref{eq:f-u} holds.
\end{enumerate}
\end{proposition}

\begin{remark}\label{lowestmodes} Note that the statement of
  Proposition \ref{evchar} does not depend on the definition of
  $\alpha_n(\lambda,k)$ for $(n+k)^2<\lambda$. In particular, Remark
  \ref{setn} follows from Proposition \ref{evchar}.
\end{remark}

\begin{proof} Let $u\in\mathcal{N}(H(k)-\lambda I)$ and write
\[ u(x,y)=\frac{1}{\sqrt{2\pi}}\sum_{n\in\Z} \hat u_n(y) e^{inx},
\qquad \hat u_n(y):=\frac{1}{\sqrt{2\pi}}\int_\T u(x,y) e^{-inx}\
dx. \] From $u\in\tilde\Upsilon$ it follows that $\hat u_n\in
H^2_{loc}(\R^2\setminus\{0\})\cap L_2(\R^2)$ and that the limits
\[ \Xi \hat u_n := -\lim_{y\rightarrow 0}
 \frac{1}{\log|y|}\ \hat u_n(y), \qquad
 \Omega \hat u_n := \lim_{y\rightarrow 0}
 \left(\hat u_n(y) +\log|y|\ \Xi \hat u_n\right), \]
exist. Moreover, $\hat u_n$ satisfies
\[ -\Delta \hat u_n = (\lambda-(n+k)^2)\hat u_n \qquad\text{in
}\R^2\setminus\{0\}. \]
It is well-known that this implies
\begin{equation*} \hat u_n (y) = \left\{
    \begin{array}{l@{\qquad\text{if}\;\;}l}
        0 & (n+k)^2<\lambda, \\
        c_n K_0(\sqrt{(n+k)^2-\lambda}\ |y|) & (n+k)^2 >\lambda.
    \end{array} \right.
\end{equation*}
for some constants $c_n\in\C$. Now (see \cite{AS})
\[ K_0(\sqrt{(n+k)^2-\lambda}\ \epsilon) + \log\epsilon \rightarrow
-2\pi\ \alpha_n(\lambda,k) \qquad (\epsilon\rightarrow 0) \]
implies that
\[ f(x) = (\Xi u)(x) = \frac{1}{\sqrt{2\pi}}\sum_{(n+k)^2>\lambda} c_n
e^{inx}, \]
\[ (\Omega u)(x) = -\sqrt{2\pi}\sum_{(n+k)^2>\lambda}\alpha_n(\lambda,k) c_n
e^{inx} =-2\pi\ (A_0(\lambda,k)f)(x). \]
In particular, we have proven that
$f\in\mathcal{D}(A_0(\lambda,k))=\mathcal{D}(A(\lambda,k))$.
Finally we conclude that $-2\pi\ A(\lambda,k)f=\Omega u -2\pi\sigma\Xi u
=0$.\\
The proof of part (2) is easier and will be omitted.
\end{proof}

We note that the operators $H(k)$ may have infinitely many
(embedded) eigenvalues.

\begin{example}\label{sigmaconst} Let $\sigma\equiv\alpha\in\R$ be
  constant. Then
\[ \sigma_p\left(H(k)\right) = \{\xi(\alpha)+(n+k)^2:\ n\in\Z\}.
\]
This follows from Proposition \ref{evchar} or directly by
separation of variables using the results from \cite{AGHH}.
\end{example}

\subsection{Proof of Proposition \ref{thmdiscpi}} Let $k$ be fixed
throughout this subsection. First we will prove that the operator
$H(k)$ has always as eigenvalue less or equal
$\xi(\tilde\sigma)+k^2$ where $\xi$ and $\tilde\sigma$ are given
by \eqref{eq:xi}, \eqref{eq:sigmaaverage}, respectively. Consider
the normalized trial function $e_0\equiv\frac{1}{\sqrt{2\pi}}\in
L_2(\T)$. Then
\[ (A(\lambda,k)e_0,e_0)=\alpha_0(\lambda,k)+\tilde\sigma. \]
This is non-positive provided $\lambda\geq
\xi(\tilde\sigma)+k^2$.\\
On the other hand, the operator $A(\lambda,k)$ is positive
definite provided $-\lambda$ is large. Since the eigenvalues of
$A(\lambda,k)$ depend continuously on $\lambda$ there is a
$\lambda_0\in(-\infty,k^2+\xi(\tilde\sigma)]$ such that
$A(\lambda_0,k)$ has eigenvalue $0$. By Proposition \ref{evchar}
this proves the first part of Proposition \ref{thmdiscpi}.\\
To prove the second part we need the following

\begin{lemma}\label{comparison} Let $\sigma_1\in L_\infty(\R)$ be
  real-valued and $2\pi$-periodic and let $H_1$, $H_1(k)$ be the
  operators corresponding to $\sigma_1$. Then $\sigma_1\leq\sigma$
  implies $H_1\leq H$, $H_1(k)\leq H(k)$.
\end{lemma}

\begin{proof} We consider only the case of the operators in
$L_2(\R^3)$. It suffices to prove that for some $a>0$ one has
$R(-a)\leq R_1(-a)$, where $R_1(z)=(H_1-zI)^{-1}$. By the identity
\eqref{eq:resdiff} and a similar identity for $R_1(-a)$ it
suffices to prove $(T(-a)+\sigma)^{-1}\leq(T(-a)+\sigma_1)^{-1}$,
which is immediate.
\end{proof}

To complete the proof of Proposition \ref{thmdiscpi} it suffices
to take $\sigma_1\equiv\text{ess-inf}\,\sigma$ and note (see Example
\ref{sigmaconst}) that the corresponding operator $H_1(k)$ has only a
finite number of eigenvalues below $k^2$. By Lemma \ref{comparison}
and the variational principle the same holds true for the operator
$H(k)$.

\subsection{Complexification}\label{sec:complex}
In this subsection we fix $k\in Q$,
$\lambda\in\R\setminus\tau(k)$ and assume in addition that $ k\neq
0$. As in \cite{FrSh} we choose $\delta\in(0,|k|)$ such that
$(n+\kappa)^2-\lambda \neq 0$ for all $n\in\Z$, $\kappa\in
[k-\delta,k+\delta]$, and note that there is a constant
$C_1=C_1(k,\lambda,\delta)>0$ such that
\begin{equation}\label{eq:analyticestimate}
    |(n+\mu)^2-\lambda| \geq C_1(1+|\im\mu|)^2,
    \qquad n\in\Z,\ \mu\in W,
\end{equation}
where we have put
\[ W:=\{\mu\in\C:\ |\re\mu-k|<\delta \}. \]
It follows that the functions $\alpha_n(\lambda,.)$, $n\in\Z$, defined
in \eqref{eq:alpha} admit an analytic extension to $W$. This allows to
define an m-sectorial operator $A(\lambda,\mu)$ for $\mu\in W$ by
\eqref{eq:defopa}. We need the following result for $\mu$ with large
imaginary part.

\begin{lemma}\label{analyticlemma} Let $k$, $\lambda$, $\delta$ be
as above. Then there exist constants $C_2=C_2(k,\lambda,\delta)$,
$\eta_0>0$ such that for all $\mu\in W$ with $|\im\mu|\geq\eta_0$ the
operator $A(\lambda,\mu)$ is boundedly invertible and
\[ \|(A(\lambda,\mu))^{-1}\| \leq \frac{C_2}{\log(1+|\im\mu|)}. \]
\end{lemma}

\begin{proof} From \eqref{eq:analyticestimate} we see that
\[ \begin{split}
  |\alpha_n(\lambda,\mu)| & \geq
    {\textstyle\frac{1}{4\pi}}\log\left|(n+\mu)^2-\lambda\right| -
    \varsigma \geq\\
    & \geq {\textstyle\frac{1}{2\pi}}\log(1+|\im\mu|) - (\varsigma
    - {\textstyle\frac{1}{4\pi}} \log C_1),
\end{split} \]
because $\varsigma>0$. We conclude that for large $|\im\mu|$ the
operator $A_0(\lambda,\mu)$ is boundedly invertible with
\[ \|(A_0(\lambda,\mu))^{-1}\| \leq \frac{C_3}{\log(1+|\im\mu|)},
\]
and we obtain the assertion by noting that
\[ (A(\lambda,\mu))^{-1}=
    \left(I+(A_0(\lambda,\mu))^{-1}\sigma\right)^{-1}(A_0(\lambda,\mu))^{-1}
\]
whenever $\|(A_0(\lambda,\mu))^{-1}\sigma\|<1$.
\end{proof}

Now we obtain easily the

\begin{proof}[Proof of Proposition \ref{thmpppi}] It suffices to
repeat the arguments in the proof of Proposition 1.10 in
\cite{Fr2}, replacing Proposition 3.5 there by our
Lemma~\ref{analyticlemma}.
\end{proof}

%%%%%%%%%%%%%%%%%%%%%%%%%%%%%%%%%%%%%%%%%%%%%%%%%%%%%%%%%%%%%%%%%%%%
%%%%%%%%%%%%%%%%%%%%%%%%%%%%%%%%%%%%%%%%%%%%%%%%%%%%%%%%%%%%%%%%%%%%

\section{The second model: an infinite family of lines}\label{sec:model2}

Now we would like to discuss a model of a ``diffraction grating''
consisting of periodically arranged wires. Our approach will be
similar to the one outlined above and we emphasize these
similarities by \emph{keeping the same notation for analogous
objects}. However, several constructions in the present case are
technically more involved and we concentrate on these in the
exposition.

\subsection{Main result for the second model}
It is convenient to denote now the coordinates in $\R^3$ by
$(x,y),\ x\in\R^2,\ y\in\R$, and to put
\[ \Gamma:= \bigcup_{n\in\Z}\{(x_1,2\pi n,0):\ x_1\in\R\}. \]
As before let $\sigma$ be a real-valued, $2\pi$-periodic function
satisfying \eqref{eq:sigma}. We will construct a  self-adjoint
operator $H$ in $L_2(\R^3)$ corresponding to the formal expression
$-\Delta+\sum_{n\in\Z}\sigma(x_1)\delta(x_2-2\pi n)\delta(y)$, see
the footnote in Subsection~\ref{sec:mainres1}. Put
$C_\epsilon:=\cup_{n\in\Z}\{(x,y)\in\R^3:\ (x_2+2\pi
n)^2+y^2\leq\epsilon^2\}$ for $\epsilon>0$. We denote by
$\Upsilon$ the class of functions $u\in L_2(\R^3)$ satisfying
\begin{equation}\label{eq:h2m2}
    u\in H^2(\R^3\setminus C_\epsilon) \qquad\text{for all
    }\epsilon>0
\end{equation}
and such that for all $n\in\Z$ the limits
\[\begin{split} \Xi_n u &:= -\lim_{(x_2,y)\rightarrow 0}
 \frac{1}{\log\sqrt{x_2^2+y^2}}\ u(.,x_2+2\pi
 n,y), \\
 \Omega_n u &:= \lim_{(x_2,y)\rightarrow 0}
 \left(u(.,x_2+2\pi n,y) +\log\sqrt{x_2^2+y^2}\ \Xi u\right)
 \end{split} \]
exist in the sense of distributions, belong to $L_2(\R)$ and satisfy
\[
\sum_{n\in\Z} \left(\|\Xi_n u\|^2+\|\Omega_n u\|^2\right) < \infty.
\]
We introduce the operators $\Xi,\Omega:L_2(\R^3)\rightarrow
\sum_{n\in\Z}^\oplus L_2(\R)$,
\[\begin{split}
(\Xi u)_n &:= \Xi_n u, \qquad  (\Omega u)_n := \Omega_n u,\\
\mathcal D(\Xi) & :=\mathcal D(\Omega):=\Upsilon.
\end{split}\]
As before (see also Subsection \ref{sec:defop1m2} below) one
constructs a self-adjoint operator $H$ in $L_2(\R^3)$ such that
\begin{equation}\label{eq:defopm2}\begin{split}
    H u & = -\Delta u \qquad\text{in }\R^3\setminus\Gamma,\\
    \mathcal{D}(H) & =\{u\in\Upsilon:\ \Delta u\in
        L_2(\R^3\setminus\Gamma),\ \Omega u - 2\pi\sigma\Xi u = 0
        \}.
    \end{split}
\end{equation}
We will again denote the standard self-adjoint realization of
$-\Delta$ in $\R^3$ by $H_0$. Our main result is

\begin{theorem}\label{thmm2} Theorems \ref{thmac},
\ref{thmscat}, \ref{thmsurf} hold also in the above situation.
\end{theorem}

\begin{remark} For simplicity we assume that our model is
  $2\pi$-periodic with respect to both $x_1$ and $x_2$. Our argument
  extends easily to the case where the periods are different. The
  assumption \eqref{eq:sigma} on $\sigma$ can be relaxed as before.
\end{remark}

\subsection{Direct integral decomposition}\label{sec:directintegralm2}
Now we write
\[ \Pi := \{ (x,y)\in\R^3:\; x\in [-\pi,\pi)^2 \} \]
and define $\tilde H^2(\Pi)$ and $\tilde\Upsilon$ in an obvious
way. As before we consider $\Xi u$, $\Omega u$ for
$u\in\tilde\Upsilon$ as functions in $L_2(\T)$.\\
For any (two-dimensional) parameter $k\in
Q:=[-\frac{1}{2},\frac{1}{2})^2$ there exists (see Subsection
\ref{sec:defoppim2}) a self-adjoint operator $H(k)$ in $L_2(\Pi)$ such
that
\begin{equation}\label{eq:defoppim2}\begin{split}
    H(k) u & =((D_x+k)^2+D_y^2) u \qquad\text{in
        }\Pi\setminus\Gamma,\\
    \mathcal{D}(H(k)) & =\{u\in\tilde\Upsilon:\ \Delta u\in
        L_2(\Pi\setminus\Gamma),\ \Omega u - 2\pi\sigma\Xi u = 0 \}.
  \end{split}
\end{equation}
Here $D_x=-i\nabla_x$, $D_y=-i\frac\partial{\partial y}$. Moreover, we
denote by $H_0(k)$ the operator $(D_x+k)^2+D_y^2$ in $L_2(\Pi)$ with
domain $\tilde H^2(\Pi)$.\\
The \emph{Gelfand transformation} $\mathcal{U}:L_2(\R^3)
\rightarrow \int_Q^\oplus L_2(\Pi)\,dk$ is in this case defined
by
\[ (\mathcal{U}u)(k,x,y) := \sum_{n\in\Z^2} e^{-i \langle k,x+2\pi
  n\rangle} u(x+2\pi n,y),
        \qquad k\in{\textstyle Q},\, (x,y)\in\Pi, \]
and realizes the unitary equivalences
\begin{equation}\label{eq:directintdecompm2}
    \mathcal{U}\, H_0\, \mathcal{U}^* = \int_Q^\oplus H_0(k)\ dk,
    \qquad \mathcal{U}\, H\, \mathcal{U}^* = \int_Q^\oplus H(k)\ dk.
\end{equation}
As before the proof of Theorem \ref{thmm2} reduces to the following

\begin{proposition}\label{thmpim2} Propositions \ref{thmacpi},
\ref{thmscpi}, \ref{thmdiscpi}, \ref{thmpppi} hold also in the
above situation.
\end{proposition}

%%%%%%%%%%%%%%%%%%%%%%%%%%%%%%%%%%%%%%%%%%%%%%%%%%%%%%%%%%%%%%%%%%%%%%%%%%%%%%%

\section{Definition of the operators in the second model}

\subsection{Definition of the operator $H$}\label{sec:defop1m2}
Let us start with some remarks concerning the definition of the
operator $H$. We consider the trace operator
$\gamma:H^2(\R^3)\rightarrow \sum_{n\in\Z}^\oplus L_2(\R)$,
\[ (\gamma u)_n(x_1) := u(x_1,2\pi n,0), \qquad x_1\in\R,\ n\in\Z.
\]
The operator $T(z)$ will in this case be an operator in
$\sum_{n\in\Z}^\oplus L_2(\R)$. We need some preparations. For
$z\in\C\setminus[0,\infty)$ we define pseudo-differential
operators $T_j(z)$, $j\in\N_0$, in $L_2(\R)$ by
\[\begin{split}
    T_0(z) &:= \frac{1}{4\pi}\log \left(D^2-z\right) -\varsigma I,\\
    T_j(z) &:= -\frac{1}{2\pi} K_0\left(2\pi j\ (D^2-z)^{1/2}\right),
    \qquad j\in\N,\\
    \mathcal{D}(T_0(z)) &:= \left\{f\in L_2(\R):\ \int_\R
    \left(\log(1+\xi^2)\right)^2 |\hat{f}(\xi)|^2\ d\xi <\infty \right\},
    \end{split}
\]
Again we choose here and in all the following the principal branches
of the logarithm and the square root on $\C\setminus(-\infty,0]$.\\
Note that $T_0(z)$ coincides with the operator $T(z)$ from
Subsection \ref{sec:defop1}. The following assertion shows in
particular that the $T_j(z)$, $j\in\N$, are bounded operators.

\begin{lemma}\label{schur} Let $z\in\C\setminus[0,\infty)$. Then
$\sum_{j\in\N}\|T_j(z)\|<\infty$.
\end{lemma}

\begin{proof} For any $\epsilon>0$ there exists a $C_\epsilon>0$
such that
\begin{equation}\label{eq:besselestimate}
    |K_0(\zeta)|\leq C_\epsilon
    \frac{e^{-\re\zeta}}{|\zeta|^{1/2}},
    \qquad |\zeta|\geq\epsilon,
\end{equation}
(see \cite{AS}). With $\epsilon:=2\pi\inf_{\xi\in\R}
|\xi^2-z|^{1/2}>0$ we find that for all $j\in\N$
\[ \|T_j(z)\| = {\textstyle\frac{1}{2\pi}} \sup_{\xi\in\R}
    |K_0(2\pi j (\xi^2-z)^{1/2})| \leq
    \tilde C_\epsilon \sup_{\xi\in\R} e^{-2\pi
    j\re\sqrt{\xi^2-z}}.
\]
Since $\re\sqrt{\xi^2-z}$ is bounded away from $0$ the assertion
follows.
\end{proof}

>From Lemma \ref{schur} and Schur's lemma one finds that the
operator $T(z)$ in $\sum_{n\in\Z}^\oplus L_2(\R)$,
\[\begin{split}
    (T(z)f)_n & := \sum_{m\in\Z} T_{|n-m|}(z)f_m, \qquad n\in\Z, \\
    \mathcal{D}(T(z)) & := \left\{f\in{\textstyle\sum_{n\in\Z}^\oplus}
    L_2(\R): \
    \sum_{n\in\Z} \int_\R \left(\log(1+\xi^2)\right)^2 |\hat
    f_n(\xi)|^2\ d\xi <\infty \right\},
    \end{split}
\]
is well-defined. Moreover, one verifies that for
$z,\zeta\in\C\setminus[0,\infty)$
\[ T(\overline z)=T(z)^*, \qquad
T(z)-T(\zeta)=(\zeta-z)\left(\gamma R_0(\zeta)\right)\left(\gamma
R_0(\overline z)\right)^*. \]
Again by \cite{Po} this implies that there exists a self-adjoint
operator $H$ in $L_2(\R^3)$ such that its resolvent
$R(z):=(H-zI)^{-1}$ is related to $R_0(z)$ by
\begin{equation}\label{eq:resdiffm2}\begin{split}
    R(z)= & R_0(z) + (\gamma R_0(\overline z))^*(T(z)+\sigma)^{-1}\gamma
    R_0(z),\\
    & z\in\C\setminus[0,\infty),\ 0\in\rho\left(T(z)+\sigma\right).
    \end{split}
\end{equation}
Here we understand $\sigma$ as an operator in $\sum_{n\in\Z}^\oplus
L_2(\R)$ acting according to $(\sigma f)_n =\sigma f_n$ for
$f=(f_n)\in\sum_{n\in\Z}^\oplus L_2(\R)$.\\
By the proof of Lemma \ref{schur} one easily finds that $T(-a)+\sigma$
is positive definite for all sufficiently large $a$, and hence $H$ is
lower semibounded.\\
Finally, one shows that $H$ satisfies \eqref{eq:defopm2}. This follows
from the identities
\[ \Xi \left(\gamma R_0(\overline z)\right)^* =
    \frac1{2\pi}I,
    \qquad \Omega \left(\gamma R_0(\overline z)\right)^* = - T(z).
    \]

\subsection{Auxiliary material}\label{sec:auxiliary}
Before we can explain the construction of the operators $H(k)$ we need
to collect some material on point interactions in a two-dimensional
strip. Note that our approach is somewhat different from the one
adopted in Section III.4 in \cite{AGHH}. Having in mind the later
application we denote the coordinates in $\R^2$ by $(x_2,y)$, the
quasi-momentum by $k_2\in Q':=[-\frac12,\frac12)$ and put
\[
\Pi':=[-\pi,\pi)\times\R.
\]
For $z\in\C\setminus[k_2^2,\infty)$ we introduce the function
\[ \begin{split}
\psi(x_2,y,z,k_2) := &
\frac12\sum_{n\in\Z}\frac{e^{inx_2}}{\sqrt{(n+k_2)^2-z}}\
e^{-\sqrt{(n+k_2)^2-z}\,|y|},\\
& (x_2,y)\in\Pi'\setminus\{(0,0)\}.
\end{split} \]
This function belongs to $L_2(\Pi')$ and is smooth away from
$(0,0)$. Moreover,

\begin{lemma}\label{psi} Let $k_2\in Q'$,
  $z\in\C$ and assume that $u\in L_2(\Pi')$ is a periodic (with
  respect to $x_2$) solution of
\begin{equation}\label{eq:twodeq}
((D_{x_2}+k_2)^2+D_y^2)u=z\, u
\qquad \text{in }\Pi'\setminus\{(0,0)\}.
\end{equation}
If $z\in[k^2,\infty)$ then $u\equiv0$, and if
  $z\in\C\setminus[k_2^2,\infty)$ then
\[ u=c\, \psi(\cdot,z,k_2), \qquad c\in\C. \]
\end{lemma}

By a \emph{periodic (with respect to $x_2$)} solution of
\eqref{eq:twodeq} we mean that the test functions in the
distributional definition of a solution are not required to vanish near
$\partial\Pi'$ but only to be periodic (with respect to $x_2$).\\
The proof of Lemma \ref{psi} follows easily by Fourier transformation
and the fact that $((D_{x_2}+k_2)^2+D_y^2-z)u$ is a distribution
supported on $\{(0,0)\}$ and hence coincides with a finite linear
combination of derivates of the $\delta$-distribution.\\
For $z\in\C\setminus[0,\infty)$ the following alternative expression
  for $\psi(\cdot,z,k)$ exists,
\begin{equation}\label{eq:psipoisson}
\psi(x_2,y,z,k_2) =
\sum_{m\in\Z} e^{-ik_2(x_2+2\pi m)}
K_0\left(\sqrt{-z}\,\sqrt{(x_2+2\pi m)^2+y^2}\right).
\end{equation}
Indeed, this follows from
\[
\frac1\pi \int K_0\left(\sqrt{-z}\,\sqrt{x_2^2+y^2}\right)\,
e^{-i\xi x_2}\, dx_2 =
\frac{e^{-\sqrt{\xi^2-z}\,|y|}}{\sqrt{\xi^2-z}},
\qquad \xi\in\R,
\]
by the Poisson summation formula.\\
Put $C_\epsilon':=\{(x_2,y)\in\R^2:\,x_2^2+y^2\leq\epsilon^2\}$. We
denote by $\tilde\Upsilon'$ the class of functions $u\in L_2(\Pi')$
satisfying
\[
u\in \tilde H^2(\Pi'\setminus C_\epsilon')
\qquad\text{for all }\epsilon>0
\]
and such that the limits
\[ \begin{split}
\Xi u &:= -\lim_{(x_2,y)\rightarrow 0}
 \frac{1}{\log\sqrt{x_2^2+y^2}}\ u(x_2,y), \\
 \Omega u &:= \lim_{(x_2,y)\rightarrow 0}
 \left(u(x_2,y) +\log\sqrt{x_2^2+y^2}\ \Xi u\right)
\end{split} \]
exist. It is not difficult to verify that
$\psi(\cdot,z,k_2)\in\tilde\Upsilon'$ with $\Xi\psi(\cdot,z,k)=1$ and
where
\[
t(z,k_2) := -\frac1{2\pi}\Omega\psi(\cdot,z,k_2),
\qquad z\in\C\setminus[k_2^2,\infty),
\]
satisfies the relation
\begin{equation}\label{eq:t}
t(z,k_2)= t(-1,0)-\frac12\sum_{n\in\Z}
\left(\frac1{\sqrt{(n+k_2)^2-z}}-\frac1{\sqrt{n^2+1}}\right).
\end{equation}
Moreover, from \eqref{eq:psipoisson} and the properties of $K_0$ (see
\cite{AS}) one deduces that for $z\in\C\setminus[0,\infty)$
\begin{equation}\label{eq:t1}
t(z,k_2)= \frac1{4\pi}\log(-z)-\varsigma
-\frac1\pi\sum_{j\in\N} \cos(2\pi jk_2) K_0(2\pi j\sqrt{-z}).
\end{equation}

We close this subsection with an estimate that will be useful in
the proof of absolute continuity. By the same arguments as in the
proof of Lemma \ref{schur} we deduce from \eqref{eq:t1} the following

\begin{lemma}\label{testimate} For any $\epsilon>0$ there is a
  constant $C=C(\epsilon)>0$ such that for all
  $z\in\C\setminus[0,\infty)$ with $\re\sqrt{-z}\geq\epsilon$ and all $k_2\in\overline{Q'}$ one has
\[
\left|t(z,k_2) -\frac1{4\pi}\log(-z)\right| \leq C.
\]
\end{lemma}

\begin{remark} Note that the subtraction of the terms $(n^2+1)^{-1/2}$
  in \eqref{eq:t} is a renormalization of the divergent sum
  $\sum((n+k_2)^2-z)^{-1/2}$. A different, but equivalent
  renormalization is chosen in Theorem III.4.8 in \cite{AGHH}.
\end{remark}

\subsection{Definition of the operators $H(k)$}\label{sec:defoppim2}
Again we denote by $\gamma$ the trace operator $\tilde
H^2(\Pi)\rightarrow L_2(\T)$. For $k\in Q$, $z\in\C\setminus
[|k|^2,\infty)$ we consider the pseudo-differential operator $T(z,k)$
  in $L_2(\T)$,
\[\begin{split}
    T(z,k) & := t(z-(D_{x_1}+k_1)^2,k_2),\\
    \mathcal{D}(T(z,k)) &:= \left\{f\in L_2(\T):\ \sum_{n\in\Z}
        \left(\log(1+n^2)\right)^2 |\hat{f}_n|^2 <\infty \right\}.
    \end{split}
\]
It follows from Lemma \ref{testimate} that $T(z,k)$ is closed and
lower semibounded and has compact resolvent. Our next goal is to show
that these operators appear as fibers in the direct integral
decomposition of the operator $T(z)$. For this purpose consider the
\emph{unitary} operator $\tilde{\mathcal U} :\sum_{n\in\Z}^\oplus
L_2(\R) \rightarrow\int_Q^\oplus L_2(\T)\, dk$, defined for smooth $f$
by
\[ (\tilde{\mathcal U} f)(k,x_1) := \sum_{m\in\Z^2} e^{-ik_1x_1-2\pi i\langle
  k,m\rangle} f_{m_2}(x_1+2\pi m_1), \qquad k\in Q,\, x_1\in\T. \]
We note that on $H^2(\R^3)$ one has the identity
\begin{equation}\label{eq:gelfandtracem2}
  \gamma\,\mathcal U = \tilde{\mathcal U}\,\gamma
\end{equation}
with an obvious meaning of the trace operator $\gamma$ on the
different sides of the equality. Moreover, it turns out that
\begin{equation}\label{eq:gelfandtm2}
  \tilde{\mathcal U}\, T(z)\, \tilde{\mathcal U}^* =
  \int_Q^\oplus T(z,k)\,dk, \qquad z\in\C\setminus[0,\infty).
\end{equation}
Similarly as in Subsection \ref{sec:defoppi} we deduce from relations
\eqref{eq:gelfandtracem2}, \eqref{eq:gelfandtm2} that there exists a
self-adjoint operator $H(k)$ in $L_2(\Pi)$ such that its resolvent
$R(z,k):=(H(k)-zI)^{-1}$ is related to $R_0(z,k):=(H_0(k)-zI)^{-1}$ by
\begin{equation}\label{eq:resdiffpim2}
    \begin{split}
    R(z,k)= & R_0(z,k)+(\gamma R_0(\overline
    z,k))^*(T(z,k)+\sigma)^{-1}\gamma R_0(z,k),\\
    & z\in\C\setminus[k^2,\infty),\
    0\in\rho\left(T(z,k)+\sigma\right),
    \end{split}
\end{equation}
and that this operator satisfies \eqref{eq:defoppim2} and
\eqref{eq:directintdecompm2}.

%%%%%%%%%%%%%%%%%%%%%%%%%%%%%%%%%%%%%%%%%%%%%%%%%%%%%%%%%%%%%%%%%%%%%%%

\section{The spectrum of the operators $H(k)$ in the second
model}

\subsection{The continuous spectrum of $H(k)$} The analogue of
Proposition \ref{thmacpi} follows exactly as before using the explicit
form of the 'unperturbed' resolvent kernel
\begin{equation*}\label{eq:freekernelm2}\begin{split}
    r_0(x,y,x',y',z,k) = & \frac{1}{8\pi^2}\sum_{n\in\Z^2}
        e^{i\langle n,x-x'\rangle}
    \frac{e^{-\sqrt{|n+k|^2-z}\,|y-y'|}}{\sqrt{|n+k|^2-z}},
        \\
    & (x,y),(x',y')\in\Pi. \end{split}
\end{equation*}
We note that the assertion of Remark \ref{spclass} remains true.\\
To establish the analogue of Proposition \ref{thmscpi} we proceed
again as before taking into account that for any $m\in\Z$, $k\in Q$,
$\lambda\in\R\setminus\tau(k)$, where now
\[
\tau(k):=\{|n+k|^2:\ n\in\Z^2\},
\]
the function $t(\cdot-(m+k_1)^2,k_2)$ has an analytic extension from
$\C_+$ to a neighbourhood of $\lambda$ in $\overline{\C_-}$. This is
easily seen from \eqref{eq:t}.

\subsection{Characterization of eigenvalues of
$H(k)$}\label{sec:opbdrym2} We turn to the point spectrum of
the operators $H(k)$ and derive an analogue of Proposition
\ref{evchar}. For $k\in Q$, $\lambda\in\R\setminus\tau(k)$ we
define
\begin{equation}\label{eq:alpham2}
    \alpha_n(\lambda,k) :=
    t(-|(n+k_1)^2+k_2^2-\lambda|+k_2^2,k_2),
    \qquad n\in\Z.
\end{equation}
In the Hilbert space $L_2(\T)$ we consider the operator
\begin{equation}\label{eq:defopam2}\begin{split}
    (A(\lambda,k)f)(x) & := \frac{1}{\sqrt{2\pi}}\sum_{n\in\Z}
    \alpha_n(\lambda,k)\hat f_n e^{inx} +\sigma(x)f(x), \qquad
    x\in\T,\\
    \mathcal{D}(A(\lambda,k)) & := \{f\in L_2(\T): \sum_{n\in\Z}
    \left(\log(1+n^2)\right)^2 |\hat{f}_n|^2 <\infty \}.
    \end{split}
\end{equation}
As for the operators $T(z,k)$ one checks that $A(\lambda,k)$ is
closed and lower semibounded and has compact resolvent. Our main tool
in the investigation of the point spectrum of $H(k)$ will be

\begin{proposition}\label{evcharm2} Let $k\in Q$ and
$\lambda\in\R\setminus\tau(k)$.
\begin{enumerate}
\item Let $u\in\mathcal{N}(H(k)-\lambda I)$ and define
\begin{equation}\label{eq:f-um2}
        f := \Xi u.
\end{equation}
Then $f\in\mathcal{N}(A(\lambda,k))$, $\hat f_n=0$ if
$(n+k_1)^2<\lambda-k_2^2$ and, moreover,
\begin{equation}\label{eq:u-fm2}
    u(x,y) = \frac{1}{\sqrt{2\pi}}\sum_{(n+k_1)^2>\lambda-k_2^2}
        \hat f_n\,e^{inx_1}\,\psi(x_2,y,(n+k_1)^2-\lambda,k_2),\
        (x,y)\in\Pi\setminus\Gamma.
\end{equation}
\item Let $f\in\mathcal{N}(A(\lambda,k))$ such that $\hat f_n=0$
if $(n+k_1)^2<\lambda-k_2^2$ and define $u$ by \eqref{eq:u-fm2}.\\
Then $u\in\mathcal{N}(H(k)-\lambda I)$ and, moreover,
\eqref{eq:f-um2} holds.
\end{enumerate}
\end{proposition}

Proceeding as in the proof of Proposition \ref{evchar} this follows
easily by Fourier transformation from Lemma \ref{psi}.\\
Again the operators $H(k)$ may have infinitely many (embedded)
eigenvalues.

\begin{example}\label{sigmaconstm2} Let $\sigma\equiv\alpha\in\R$ be
constant. It follows from \eqref{eq:t} and Lemma \ref{testimate} that
the function $t(\cdot,k_2)$ is continuous and decreasing on
$(-\infty,k_2^2)$ with
\[
\lim_{\lambda\to k_2^2-} t(\lambda,k_2)=-\infty, \qquad
\lim_{\lambda\to-\infty} t(\lambda,k_2)=\infty.
\]
Hence there exists a unique $\lambda(\alpha,k_2)\in(-\infty,k_2^2)$
such that
\[
t(\lambda(\alpha,k_2),k_2)+\alpha=0
\]
>From Proposition \ref{evcharm2} it follows that
\[
\sigma_p\left(H(k)\right)=\{\lambda(\alpha,k_2)+(n+k_1)^2:\ n\in\Z\}.
\]
This result (in equivalent notation) may also be deduced from the
two-dimensional result in \cite{AGHH} by separation of variables.
\end{example}

For the proof of the analogue of Proposition \ref{thmdiscpi} we
proceed exactly as before. Using the trial function $e_0$ we find that
$H(k)$ has an eigenvalue less or equal
$\lambda(\tilde\sigma,k_2)+k_1^2$, where $\tilde\sigma$ is given by
\eqref{eq:sigmaaverage} and $\lambda(\cdot,k_2)$ was constructed in
Example \ref{sigmaconstm2}. Since $H_1(k)$, the operator corresponding
to $\sigma_1\equiv\text{ess-inf}\,\sigma$, has only finitely many
eigenvalues an analogue of Lemma \ref{comparison} finishes the proof
of the analogue of Proposition \ref{thmdiscpi}.
\bigskip

Finally, we turn to the proof of the analogue of Proposition
\ref{thmpppi}. Again we will construct an analytic extension of the
operators $A(\lambda,k_1,k_2)$ with respect to the variable $k_1$. The
new ingredient here is that we replace the 'complicated' symbol
$\alpha_n(\lambda,k)$ by the explicit
$(4\pi)^{-1}\log((n+k_1)^2-\lambda)$ (which already appeared in the
first part of our analysis). This requires careful estimates which are
uniform in $n$ and in the (complex) parameter $k_1$.\\
Now we describe the details of our construction. Similarly as in
Subsection \ref{sec:complex} we fix $k\in Q$,
$\lambda\in\R\setminus\tau(k)$ and assume that $k_1\neq 0$. We choose
$\delta\in(0,|k_1|)$ satisfying $(n+\mu)^2+k_2^2-\lambda\neq 0$ for
all $n\in\Z$, $\kappa\in[k_1-\delta,k_1+\delta]$, and put
\[
W:=\{\mu\in\C:\ |\re\mu-k_1|<\delta\}.
\]
For $n\in\Z$ we consider the functions
\begin{equation*}
  s_n(\mu):= \left\{
  \begin{array}{l@{\qquad\text{if}\;\;}l}
    -(n+\mu)^2+\lambda, & (n+k_1)^2+k_2^2-\lambda>0,\\
    (n+\mu)^2-\lambda+2k_2^2, & (n+k_1)^2+k_2^2-\lambda<0.
  \end{array} \right.
\end{equation*}
We need the technical

\begin{lemma}\label{analyticre} There are constants
  $\eta_0$, $\epsilon>0$ such that for all $n\in\Z$, $\mu\in W$ with
  $|\im\mu|\geq\eta_0$ one has $s_n(\mu)\in\C\setminus[0,\infty)$
  and
\[
\re\sqrt{-s_n(\mu)}\geq \epsilon(1+|n|),
\qquad n\in\Z,\ \mu\in W,\ |\im\mu|\geq\eta_0.
\]
\end{lemma}

\begin{proof} We will write $\mu=\kappa+i\eta$ with
  $\kappa,\eta\in\R$. Introduce
\[
J:=\{n\in\Z:\ (n+\kappa)^2>\lambda \quad \forall |\kappa-k_1|\leq\delta\}.
\]
The set $\Z\setminus J$ is finite and for $n$ from that set the
assertion is readily verified. Hence we will now consider only $n\in
J$. In particular, we have
$s_n(\mu)=-(n+\mu)^2+\lambda$. Moreover, we will use the
elementary estimate
\begin{equation}\label{eq:analyticestimate1}
    |(n+\mu)^2-\lambda| \geq c_1 \min\{(1+|n|)^2,(1+|\eta|)^2\} ,
    \qquad n\in J,\ \mu\in W.
\end{equation}
First, we assume that
\begin{equation}\label{eq:analyticre1}
\re s_n(\mu)=\eta^2-(n+\kappa)^2+\lambda\leq 0,
\end{equation}
and consequently
\[
\re\sqrt{-s_n(\mu)}=
\sqrt{|(n+\mu)^2-\lambda|}\ \cos\left(\frac12\arctan
\frac{2\eta(n+\kappa)}{(n+\kappa)^2-\eta^2-\lambda}\right)
\]
(where $\arctan(\pm\infty)=\mp\frac\pi2$). Noting that the cosine
factor is bounded away from $0$ the assertion follows immediately
from
\eqref{eq:analyticestimate1}.\\
Now we assume that the opposite inequality in \eqref{eq:analyticre1}
holds and consequently
\[
\re\sqrt{-s_n(\mu)}=
\sqrt{|(n+\mu)^2-\lambda|}\ \sin\left(\frac12\arctan
\frac{2|\eta(n+\kappa)|}{\eta^2-(n+\kappa)^2+\lambda}\right).
\]
If, say, $2|\eta(n+\kappa)|\geq\eta^2-(n+\kappa)^2+\lambda$ then
the sine factor is bounded away from $0$ and the assertion follows
again by \eqref{eq:analyticestimate1}. Otherwise if
\mbox{$2|\eta(n+\kappa)|$} $<\eta^2-(n+\kappa)^2+\lambda$ we note
that there is a constant $c_2>0$ such that $\sin(\frac12\arctan
x)\geq c_2x$ for all $0\leq x\leq 1$. Hence
\[
\re\sqrt{-s_n(\mu)} \geq 2c_2 \sqrt{|(n+\mu)^2-\lambda|}\
\frac{|\eta(n+\kappa)|}{\eta^2-(n+\kappa)^2+\lambda}.
\]
Using now \eqref{eq:analyticestimate1} and the estimate
\[
\frac{|\eta(n+\kappa)|}{\eta^2-(n+\kappa)^2+\lambda} \geq
\frac{|n+\kappa|}{|\eta|}\geq c_3 \frac{1+|n|}{|\eta|}
\]
we obtain the assertion.
\end{proof}

Since the $s_n$ assume values in $\C\setminus[k_2^2,\infty)$, it
  follows from our discussion in Subsection \ref{sec:auxiliary} that
  $\alpha_n(\lambda,\mu,k_2):=t_n(s_n(\mu),k_2)$ defines an analytic
  function of $\mu\in W$. This allows to define an m-sectorial
  operator $A(\lambda,\mu)$ for $\mu\in W$ by \eqref{eq:defopam2}. The
  analogue of Proposition \ref{thmpppi} will follow as before if we
  can establish

\begin{lemma}\label{analyticlemmam2} Lemma \ref{analyticlemma} holds
  also in the above situation.
\end{lemma}

\begin{proof} By Lemmas \ref{testimate} and \ref{analyticre} we have
\[
\left|\alpha_n(\lambda,\mu,k_2)-
\frac1{4\pi}\log(-s_n(\mu))\right| \leq C,
\qquad n\in\Z,\ \mu\in W,\ |\im\mu|\geq\eta_0.
\]
Using this one can proceed as in the proof of Lemma
\ref{analyticlemma}.
\end{proof}

\begin{remark} Assume $\lambda<k_1^2$ for simplicity. A closer look at
  the proof of Lemma \ref{testimate} yields that we have
\[
\left|\alpha_n(\lambda,k_1,k_2)-
\frac1{4\pi}\log\left((n+k_1)^2-\lambda\right)+\varsigma\right|
\leq Ce^{-c|n|}, \qquad n\in\Z,
\]
(even uniformly in $k_1\in W$). This reflects the physical fact
that the interaction of the wires, being of a tunneling nature,
decreases exponentially fast with their distance $2\pi|n|$.
\end{remark}

\subsection*{Acknowledgments}

The research was supported in part by ASCR and its Grant Agency
within the projects IRP AV0Z10480505 and A100480501. The second
author acknowledges gratefully a partial support through the ESF
SPECT programme.

%%%%%%%%%%%%%%%%%%%%%%%%%%%%%%%%%%%%%%%%%%%%%%%%%%%%%%%%%%%%%%%%%%%%%%%%%%%%%%%
%%%%%%%%%%%%%%%%%%%%%%%%%%%%%%%%%%%%%%%%%%%%%%%%%%%%%%%%%%%%%%%%%%%%%%%%%%%%%%%

\bibliographystyle{amsalpha}

\end{document}